\documentclass{amsart}

\usepackage{amsmath,amssymb,graphicx}

\newtheorem{theorem}{Theorem}
\newcommand{\bt}{\begin{theorem}}
\newcommand{\et}{\end{theorem}}
\newtheorem{lemma}{Lemma}
\newcommand{\bl}{\begin{lemma}}
\newcommand{\el}{\end{lemma}}
\newtheorem{corollary}{Corollary}
\newcommand{\bc}{\begin{corollary}}
\newcommand{\ec}{\end{corollary}}
\newcommand{\beq}{\begin{equation}}
\newcommand{\eeq}{\end{equation}}
\newcommand{\benum}{\begin{enumerate}}
\newcommand{\eenum}{\end{enumerate}}
\newcommand{\N}{\ensuremath{ \mathbf N }}
\newcommand{\Z}{\ensuremath{\mathbf Z}}

\newcommand{\mcc}{\ensuremath{ \mathcal C}}
\newcommand{\mcd}{\ensuremath{ \mathcal D}}

\DeclareMathOperator{\qand}{\quad\text{and}\quad}
\DeclareMathOperator{\qqand}{\qquad\text{and}\qquad}

\title[Asymptotic approximate groups]{Every finite set of integers is an asymptotic approximate group}
\author{Melvyn B. Nathanson}
\address{Department of Mathematics\\Lehman College (CUNY)\\Bronx, NY 10468} \email{melvyn.nathanson@lehman.cuny.edu}

\thanks{Supported in part by a grant from the PSC-CUNY Research Award Program.}

\date{\today}

\begin{document}

\begin{abstract}
A  set $A$  is an \emph{$(r,\ell)$-approximate group} in the additive abelian group $G$ 
if $A$ is a nonempty subset of  $G$ and there exists a subset $X$ 
of $G$ such that $|X| \leq \ell$ and $rA \subseteq X+A$.
The set $A$ is an \emph{asymptotic $(r,\ell)$-approximate group} 
if the sumset $hA$ is an $(r,\ell)$-approximate group for all sufficiently large integers $h$.  
It is proved that every finite set of integers is an asymptotic  $(r,r+1)$-approximate group 
for every  integer $r \geq 2$.
\end{abstract}

\maketitle

\section{Asymptotic approximate groups}
Let $A$ be a nonempty subset of an additive abelian group $G$.   
We do not assume that $A$ is finite, nor that $A$ contains the identity, 
nor that $A$ is symmetric in the sense that $-A=A$, where $-A = \{ -a: a \in  A\}$.  
Let $x \in G$.  
The \index{translate}\emph{translate of $A$ by $x$} 
is the set $x+A = \{x+a:a\in A\}$. 
If $A_1,A_2,\ldots, A_h$ are nonempty subsets of $G$, 
then their \emph{sumset} is 
\[
A_1 + A_2 + \cdots + A_h = \{a_1 + a_2 + \cdots + a_h : 
a_i \in A_i \text{ for } i=1,2,\ldots, h\}.
\]
If $A_i = A$ for all $i=1,\ldots, h$, then we denote the sumset 
$A_1 + A_2 + \cdots  + A_h$ by $hA$.

Let $r$ and $\ell$ be positive integers.  
The  set $A$  is an \emph{$(r,\ell)$-approximate group} 
if there exists a subset $X$ of $G$ such that 
\beq              \label{AsympApproxInt:additive}
|X| \leq \ell \qqand rA \subseteq X+A.
\eeq
In this paper we consider only abelian groups because  
the subject is approximate groups of integers.  
Our definition of approximate group is more general 
than the original definition,  which is due to Tao~\cite{breu13,gree12,tao07} 
and which has been extensively investigated.

We begin with some simple observations.  
Let $A$ and $X$ be subsets of $G$ that satisfy~\eqref{AsympApproxInt:additive}.  
Let $y \in G$, and let 
\[
A' = y + A \qand X' = (r-1)y+X.
\]
Then
\[
|X'| = |X| \leq \ell
\]
and
\begin{align*}
rA' & = r(y+A) = ry + rA \\
& \subseteq ry+X+A \\
& = (r-1)y+X + (y+A) \\
& = X' + A'.
\end{align*}
Thus, every translate of an $(r,\ell)$-approximate group 
is an $(r,\ell)$-approximate group.

Let $A$ and $X$ be subsets of $G$ that satisfy~\eqref{AsympApproxInt:additive}, 
and let $a_0 \in A$.  For  $r' \in \{1,2,\ldots, r-1\}$, let  
\[
X_{r'} = -(r-r')a_0 + X.
\]
We have 
\[
|X_{r'}| = |X| \leq \ell
\]
and
\[
(r-r')a_0 + r'A \subseteq rA  \subseteq X+A.
\]
It follows that 
\[
r'A  \subseteq (-(r-r')a_0 + X) +A = X_{r'} +A
\]
and so $A$ is also a $(r',\ell)$-approximate group.

An $(r,\ell)$-approximate group is not necessarily an $(r',\ell)$-approximate  
group for $r' > r$.
If $A$ and $X$ are finite subsets of $G$ that satisfy~\eqref{AsympApproxInt:additive}, 
then 
\[
|rA| \leq |X+A| \leq \ell |A|.
\]
Equivalently, if $|rA| > \ell |A|$, then $A$ is not an $(r,\ell)$-approximate group.
For example, $\{0,1,2,3\}$ is a $(2,2)$-approximate group of integers, but not a 
$(3,2)$-approximate group.

In this paper we study asymptotic approximate groups.
The  set $A$  is an \emph{asymptotic  $(r,\ell)$-approximate group} 
if $hA$ is an $(r,\ell)$-asymptotic group for every sufficiently large integer $h$.
This means that, for every integer $h \geq h_1(A)$, there exists a subset 
$X_h$ of $G$ such that 
\beq              \label{AsympApproxInt:asymp}
|X_h| \leq \ell \qqand rh A \subseteq X_h+ hA.  
\eeq
The purpose of this paper is to prove that every nonempty 
finite set of integers is an asymptotic 
$(r,r+1)$-approximate group for every $r \geq 2$.

\section{Approximate groups of integers}
Let \Z\ and  \N\ denote, respectively, the additive group of integers and the 
additive semigroup of positive integers.
For $A \subseteq \Z$ and $d \in \Z$, we define the \emph{dilation} 
\[
d \ast A = \{da:a \in A\}.
\]
For  $u,v \in \Z$, the set 
\[
[u,v] = \{n\in \Z: u \leq n \leq v \}
\]
is an  \emph{interval of integers} .  
%Let $A$ be a set of integers, and let $u,v \in \Z$.  
The interval of integers $[u,v]$ is a \emph{maximal subinterval} of $A$ if $u-1 \notin A$, 
$v+1 \notin A$, and $[u,v] \subseteq A$.

We recall the following well-known result.

\bt                      \label{AsympApproxInt:theorem:sumsetIneq}
Let $A$ be a nonempty finite set of integers.  
For every positive integer $r$,   
\[
|rA| \geq  r|A| - r+1 = |A| + (r-1)(|A| - 1)  
\]
and
\[
|rA| =  r|A| - r+1  
\]
if and only if $A$ is an arithmetic progression.  
\et

\begin{proof}
See~\cite[Theorem 1.3]{nath96bb}.
\end{proof}

\bt
A nonempty finite set $A$ of integers is an $(r,1)$-approximate group 
for some integer $r \geq 2$ if and only if $|A| = 1$ 
if and only if $A$ is an $(r,1)$-approximate group for all $r \geq 2$.  
\et

\begin{proof}
If $A$ is an $(r,1)$-approximate group, then there is 
a set $X$ with $|X| = 1$ such that $rA \subseteq X+A$, and so
\[
 |rA| \leq |X+A| \leq |X| |A| = |A|.
\]
Theorem~\ref{AsympApproxInt:theorem:sumsetIneq} implies that 
that $|rA| > |A|$ if $r \geq 2$ and $|A| \geq 2$.  
It follows that if $r \geq 2$ and $A$ is an $(r,1)$-approximate group, then $|A| = 1$.

If $|A| = 1$ and $A = \{a_0\}$, then 
\[
rA = \{ra_0\} = \{ (r-1)a_0\} + \{a_0\} = X+A
\]  
where $X =  \{ (r-1)a_0\}$ and $|X| = 1$.  
Thus, if  $|A| = 1$, 
then $A$ is an $(r,1)$-approximate group for all $r \geq 2$.  
This completes the proof.  
\end{proof}

\bt
Let $A$ be a finite set of integers with  $|A| \geq 2$.  
The set $A$  is an $(r,2)$-approximate group for some integer $r \geq 2$ 
if and only if $r=2$ and $A$ is an arithmetic progression, 
or $r=3$ and $|A| = 2$. 
\et

\begin{proof}
Let $|A| = k \geq 2$.  
If $r \geq 2$ and  $A$ is an $(r,2)$-approximate group, then there is a set $X$ 
with $|X| = 2$ such that 
\[
rA \subseteq X+A.
\]
By Theorem~\ref{AsympApproxInt:theorem:sumsetIneq},  
\[
rk-r+1 \leq |rA| \leq |X+A| \leq |X| |A| = 2k
\]
and so 
\[
(r-2)(k-1) \leq 1,
\]
If $r \geq 3$, then $r=3$ and $k=2$.
Otherwise, $r=2$ and $A$ is a $(2,2)$-approximate group.

Note that if  $r=3$ and if $k=2$ and $A = \{a_0,a_1\}$ with $a_0 < a_1$, then 
\[
3A = \{3a_0,2a_0 +a_1,a_0+2a_1,3a_1\}
= \{2a_0, 2a_1 \}+\{ a_0,a_1\} = X+A
\]
where $X = \{2a_0, 2a_1\}$.
Thus, every set $A$ with $|A| = 2$ is a $(3,2)$-approximate group. 

Suppose that $A$ is a $(2,2)$-approximate group, 
and that $2A \subseteq X+A$, where $|X| = 2$.  
If $A$ is not an arithmetic progression, then, 
by Theorem~\ref{AsympApproxInt:theorem:sumsetIneq}, 
\[
2k \leq |2A| \leq |X| |A| = 2k
\]
and so $|2A| = |X+A| = 2k$.
Let $X = \{x_0,x_1\}$, where $x_0 < x_1$.
Then
\[
2A = X+A = (x_0 + A) + (x_1 + A)
\]
and
\[
(x_0 + A) \cap (x_1 + A) = \emptyset.  
\]
Let $a_0 = \min(A)$ and $a_{k-1} = \max(A)$.
Then
\[
2a_0 = \min(2A) = \min(X+A) = x_0 + a_0
\]
and so $a_0 = x_0$.  Similarly,
\[
2a_{k-1} = \max(2A) = \max(X+A) = x_1+a_{k-1}
\]
and so $x_1 = a_{k-1}$.
It follows that 
\[
a_0+a_{k-1} \in (x_0 + A) \cap (x_1 + A) = \emptyset
\]
which is absurd.  Therefore, if the finite set $A$ is a $(2,2)$-approximate group, 
then $A$ is an arithmetic progression.  

We shall prove that every finite arithmetic progression is a $(2,2)$-approximate group.  
If 
\[
A = \{a_0 +id:i=0,1,\ldots, k-1\}
\]
then 
\begin{align*}
2A & = \{ 2a_0 +id:i=0,1,\ldots, 2k-2\} \\
& = \{ 2a_0 +id:i=0,1,\ldots, k-1 \} \cup \{ 2a_0 + id:i= k,\ldots, 2k-1 \} \\
& \subseteq \left(   \{a_0 \} + A \right) \cup \left( \{a_0 +kd\} + A  \right) \\
& = X+A
\end{align*}
where
\[
X = \{a_0+d, a_0+kd\}
\]
and so $A$ is a 2-approximate group.  
This completes the proof.  
\end{proof}

\bt   \label{AsympApproxInt:theorem:3-linear}
Let $u_0,u,v,v_0$ be integers such that 
\[
u_0 \leq u \leq v \leq v_0
\]
and let $A$ be a set of integers such that 
\[
\{  u_0, v_0 \} \cup [u,v] \subseteq A \subseteq [u_0,v_0].
\]
Let $r \geq 2$.  If the integer $\ell$ satisfies 
\beq         \label{AsympApproxInt:ell-inequality}
\ell \geq \frac{r(v_0 - u_0) + 1}{ v-u+1 }
\eeq
and if
\[
X_{\ell} = \left\{ ru_0 - u+i(v-u+1):i=0,1,2,\ldots, \ell-1 \right\}
\]
then 
\[
rA \subseteq X_{\ell} + A
\]
and $A$ is an $(r,\ell)$-approximate group.  
\et

\begin{proof}
For $i= 0,1,\ldots, \ell - 1$ we have
\[
ru_0 + v - u +i(v-u+1) = ru_0+ (i+1)(v-u+1) - 1
\]
and so
\begin{align*} 
X_{\ell} + [u,v] 
& = \bigcup_{i=0}^{\ell-1} (ru_0 - u +i(v-u+1)) + [u,v] \\
& = \bigcup_{i=0}^{\ell-1}  [ru_0+i(v-u+1)   ,   ru_0  + v - u +i(v-u+1) ] \\
& = \bigcup_{i=0}^{\ell-1}  [ru_0+i(v-u+1)   , ru_0+ (i+1)(v-u+1) - 1 ] \\
& =  [ ru_0 , ru_0+ \ell (v-u+1) -1].
\end{align*}
Because 
\[
X_{\ell} + [u,v] \subseteq X_{\ell} +A
\]
it follows that
\[
rA  \subseteq [ru_0, rv_0]  \subseteq X_{\ell} + [u,v]
\]
if 
\[
ru_0+ \ell (v-u+1) -1 \geq rv_0.
\]
This completes the proof.  
\end{proof}

\section{Asymptotic approximate groups of integers}
The following is a fundamental result of additive number 
theory~\cite{mist-pand14,nath72e,wu-chen-chen11,yang11,yang-chen15}.

\bt     \label{AsympApproxInt:FTANT}
If  $A^{(N)}$ is a nonempty finite set of nonnegative integers 
with 
\[
\min\left( A^{(N)} \right)= 0, \quad \max \left( A^{(N)} \right) = a^*,
\qand \gcd\left( A^{(N)} \right)= 1
\]
then there is a positive integer $h_0\left( A^{(N)} \right)$ 
and there are nonnegative integers $C$ and $D$ and finite sets 
$\mcc \subseteq [0,C-2]$ and $\mcd \subseteq [0,C-2]$
such that $[C,ha^* - D]$ is a maximal subinterval of $hA^{(N)}$ and 
\[
 \{ 0, ha^* \} \cup [ C,ha^* - D] \subseteq hA^{(N)} = \mcc \cup [ C,ha^* - D] \cup (ha^* - \mcd) 
 \subseteq [0, ha^*]
\]
for all $h \geq h_0 \left( A^{(N)} \right)$.  
\et

Note that $0 \in \mcc$ if $\mcc \neq \emptyset$ and   
$0 \in \mcd$ if $\mcd \neq \emptyset$.    

\begin{proof}
See~\cite[Theorem 1.1]{nath96bb}.
\end{proof}

\bt                       \label{AsympApproxInt:theorem:FiniteInteger}
Every nonempty finite set of integers is an asymptotic $(r,r+1)$-approximate group
for every $r \geq 2$.
\et

\begin{proof}
We begin with a simple observation.  
Let $A^{(N)}$ be an $(r,\ell)$-approximate group, and let $X'$ be a set of integers 
such that $|X'| = \ell$ and 
\[
rA^{(N)} \subseteq X'+A^{(N)}.
\]
For $d \in \N$ and $a_0 \in \Z$, let
\[
A= d\ast A^{(N)} + a_0
\]
and
\[
X = d\ast X' + (r-1)a_0.
\]
Then $|X| = |X'| = \ell$, and 
\begin{align*}
rA & = r(d\ast A^{(N)} + a_0) = d\ast rA^{(N)} + ra_0 \\
& \subseteq d\ast (X'+ A^{(N)}) + ra_0 
= d\ast X' +  d\ast A^{(N)} + ra_0  \\
& = d\ast X' + (r-1)a_0 + d\ast A^{(N)} + a_0  \\ 
& = X + A.
\end{align*}
It follows that $A$ is an $(r,\ell)$-approximate group.

Let $A$ be a nonempty finite set of integers, and let $a_0 = \min(A)$ and $d = \gcd(A-a_0)$.  
The finite set 
\[
A^{(N)} = \left\{ \frac{a-a_0}{d} : a \in A \right\}
\]
satisfies $\min(A^{(N)}) = 0$,  $\gcd(A^{(N)}) = 1$, 
and 
\[
hA= d\ast hA^{(N)} + ha_0 
\]
for all $h \in \N$.  
By the previous remark, it suffices to prove that $A^{(N)}$ 
is an asymptotic $(r,r+1)$-approximate group.  

Let $a^*$, $C$, $D$, and $h_0(A^{(N)})$ be the integers defined by 
Theorem~\ref{AsympApproxInt:FTANT}, 
and let 
\[
h_1(A^{(N)}) = \max\left( h_0(A^{(N)}), \frac{ (r+1)(C+D)-r}{a^*} \right).
\]
Let   $h \geq h_1(A^{(N)})$.  
Applying Theorem~\ref{AsympApproxInt:theorem:3-linear} 
to the set $hA^{(N)}$ with $u_0 = 0$, $u=C$, 
$v = ha^*-D$, and $v_0 = ha^*$,
we satisfy inequality~\eqref{AsympApproxInt:ell-inequality} 
with $\ell = r+1$ and obtain a set 
\[
X_{h,r+1}  = \left\{-C+i(ha^*- C - D+1):i=0,1,2,\ldots, r \right\}
\]
such that $|X_{h,r+1}| = r+1$ and 
\[
rhA^{(N)} \subseteq  X_{h,r+1} + hA^{(N)}.
\]
Thus,  $A^{(N)}$ is an asymptotic $(r,r+1)$-approximate group.  
This completes the proof.  
\end{proof}

\emph{Remark.}  Using results about convex polytopes, 
Nathanson~\cite{nath15}  proved that every finite set of lattice points 
is an asymptotic approximate group, and, more generally, 
that every finite subset of an abelian group 
is an asymptotic approximate group.

\def\cprime{$'$} \def\cprime{$'$} \def\cprime{$'$} \def\cprime{$'$}
\providecommand{\bysame}{\leavevmode\hbox to3em{\hrulefill}\thinspace}
\providecommand{\MR}{\relax\ifhmode\unskip\space\fi MR }
\providecommand{\MRhref}[2]{%
  \href{http://www.ams.org/mathscinet-getitem?mr=#1}{#2}
}
\providecommand{\href}[2]{#2}


\begin{thebibliography}{1}

\bibitem{breu13}
Emmanuel Breuillard, \emph{A brief introduction to approximate groups}, Thin
  groups and superstrong approximation, Math. Sci. Res. Inst. Publ., vol.~61,
  Cambridge Univ. Press, Cambridge, 2014, pp.~23--50.

\bibitem{gree12}
Ben Green, \emph{What is{$\ldots$} an approximate group?}, Notices Amer. Math.
  Soc. \textbf{59} (2012), no.~5, 655--656.

\bibitem{mist-pand14}
Raj~Kumar Mistri and Ram~Krishna Pandey, \emph{A generalization of sumsets of
  set of integers}, J. Number Theory \textbf{143} (2014), 334--356.

\bibitem{nath72e}
Melvyn~B. Nathanson, \emph{Sums of finite sets of integers}, Amer. Math. Monthly
  \textbf{79} (1972), 1010--1012.

\bibitem{nath96bb}
\bysame, \emph{{Additive Number Theory: Inverse Problems and the Geometry of
  Sumsets}}, Graduate Texts in Mathematics, vol. 165, Springer-Verlag, New
  York, 1996.
  
  

\bibitem{nath15}
\bysame, \emph{Every finite set of integers is an asymptotic approximate group},
arXiv:1512.03130, 2015.

\bibitem{tao07}
Terence Tao, \emph{Product set estimates for non-commutative groups},
  Combinatorica \textbf{28} (2008), no.~5, 547--594.

\bibitem{wu-chen-chen11}
Jian-Dong Wu, Feng-Juan Chen, and Yong-Gao Chen, \emph{On the structure of the
  sumsets}, Discrete Math. \textbf{311} (2011), no.~6, 408--412.

\bibitem{yang11}
Quan-Hui Yang, \emph{Another proof of {N}athanson's theorems}, J. Integer Seq.
  \textbf{14} (2011), no.~8, Article 11.8.4, 5.

\bibitem{yang-chen15}
Quan-Hui Yang and Yong-Gao Chen, \emph{On the cardinality of general {$h$}-fold
  sumsets}, European J. Combin. \textbf{47} (2015), 103--114.

\end{thebibliography}
\end{document}